\documentclass[12pt,draftcls,onecolumn]{IEEEtran}
\usepackage{amsfonts,amssymb,epsfig, float}
\usepackage{times}

\newcommand{\ssm}{\scriptscriptstyle}
\newcommand{\implies}{\Rightarrow}
\newcommand{\eps}{\varepsilon}

\newcommand{\rref}[1]{(\ref{#1})}

\newtheorem{theorem}{Theorem}
\newtheorem{itlemma}{Lemma}[section] 
\newtheorem{itproposition}[itlemma]{Proposition}
\newtheorem{itcorollary}[itlemma]{Corollary}
\newtheorem{itdefinition}[itlemma]{Definition}
\newtheorem{itexample}[itlemma]{Example}

\newenvironment{lemma}{\begin{itlemma}\rm}{\end{itlemma}} 
\newtheorem{remark}[theorem]{Remark}
\newenvironment{corollary}{\begin{itcorollary}\rm}{\end{itcorollary}}
\newenvironment{proposition}{\begin{itproposition}\rm}{\end{itproposition}}
\newenvironment{definition}{\begin{itdefinition}\rm}{\end{itdefinition}}
\newenvironment{example}{\begin{itexample}\rm}{\end{itexample}}

\newcommand{\text}[1]{\hbox{\rm \ #1\ \/}}
\newcommand{\uc}{{\bf u}}
\newcommand{\vc}{{\bf v}}
\newcommand{\bl}[1]{\begin{lemma}\label{#1}}
\newcommand{\br}[1]{\begin{remark}\label{#1}}
\newcommand{\bt}[1]{\begin{theoreremar}\label{#1}}
\newcommand{\bd}[1]{\begin{definition}label{#1}}
\newcommand{\bp}[1]{\begin{proposition}\label{#1}}
\newcommand{\bc}[1]{\begin{corollary}\label{#1}}
\newcommand{\bfact}[1]{\begin{fact}\label{#1}}
\newcommand{\bex}[1]{\begin{example}\label{#1}}
\newcommand{\bem}[1]{\begin{example}\label{#1}}  
\newcommand{\ec}{\end{corollary}}
\newcommand{\eex}{\end{example}}
\newcommand{\eem}{\end{example}}
\newcommand{\el}{\end{lemma}}
\newcommand{\er}{\end{remark}}
\newcommand{\et}{\end{theorem}}
\newcommand{\ed}{\end{definition}}
\newcommand{\ep}{\end{proposition}}
\newcommand{\epr}{\end{proof}}
\newcommand{\bpr}{\begin{proof}}

\newcommand{\beq}{\begin{eqnarray}}
\newcommand{\eeq}{\end{eqnarray}}
\newcommand{\beqn}{\begin{eqnarray*}}
\newcommand{\eeqn}{\end{eqnarray*}}
\newcommand{\bi}{\begin{itemize}}
\newcommand{\ei}{\end{itemize}}
\newcommand{\ben}{\begin{enumerate}}
\newcommand{\een}{\end{enumerate}}

\newcommand{\R}{{\mathbb R}}  
\newcommand{\N}{{\mathbb N}}  

\title{A Small-Gain Theorem for Monotone Systems with Multi-Valued
Input-State Characteristics
}
\author{\vspace{.2em}Patrick De Leenheer \footnote{De Leenheer (Corresponding Author):   Tel.: +1 352 392 0281 ext.
240; Fax: +1 352 392 8357; Department of Mathematics; University
of Florida; 411 Little Hall; PO Box 118105; Gainesville, FL
32611--8105 USA; deleenhe@math.ufl.edu. Supported by NSF/DMS Grant
0500861.}
 \and Michael
Malisoff \footnote{Malisoff: Department of Mathematics; Louisiana
State University; 304 Lockett Hall; Baton Rouge, LA 70803-4918
USA; malisoff@lsu.edu. Supported by NSF/DMS Grant 0424011.}}
\begin{document}
\maketitle
\begin{abstract}
We provide a new global small-gain theorem for feedback
interconnections of monotone input-output systems with
multi-valued input-state characteristics. This extends a recent
small-gain theorem of Angeli and Sontag for monotone systems with
singleton-valued characteristics.  We prove our theorem using
Thieme's convergence theory for asymptotically autonomous systems.
An illustrative example is also provided.\smallskip\medskip


\end{abstract}

\section{Introduction}
The recent extension \cite{AS03} of the theory of monotone
dynamical systems to monotone input-output (i/o) systems has
proven to be very useful in analyzing the global behavior of many
important dynamics; see for example \cite{AS03, AS04, AS04a, SCL04, POSTA03, MBE05}, and
see Section \ref{sec2} below for the relevant definitions.  (See
also \cite{S95} for a detailed account of monotone dynamical
systems.) Of particular interest in this literature are feedback
interconnections of subsystems--or ``modules''--that are monotone
and that possess a unique globally asymptotically stable
equilibrium, obviously depending on the particular (constant)
input applied. This has lead to the introduction of the notion of
{\em input-state (i/s) characteristics}, which are maps assigning
to each constant input value the particular equilibrium point to
which solutions converge. In many applications, this {}assignment
is exactly the type of quantitative information that is available
from experiments (such as gene expression levels, for instance).
Monotonicity,{} on the other hand, may be considered as a
qualitative or structural property of an i/o system; see the
graphical tests for monotonicity in \cite{AS04} for example. These
two ingredients, monotonicity of the subsystems and existence of
characteristics, are key to proving the small-gain theorems in
\cite{AS03, AS04, AS04a, SCL04}.
(For small-gain theorems for nonlinear but not necessarily
monotone systems, see \cite{zhong-ping}.)

In practice however, many monotone i/o systems subject to constant
inputs possess {\em several} equilibria and all solutions converge
to one of them, although distinct solutions may converge to
distinct equilibria. Such systems are sometimes called {\it
multi-stable}. In fact, since monotone i/o systems subject to
constant inputs are monotone dynamical systems,
 this type of
global behavior is to be expected (see \cite{S95}). This suggests
that the notion of an i/s characteristic ought to  be generalized
to a {\it multi-valued} map which assigns to each constant input
value the set of all possible equilibria to which solutions
converge.

This naturally leads to  the question of  whether the known
small-gain theorem for monotone systems in \cite{AS03} remains
valid if instead of the original notion of i/s characteristics,
one assumes the existence of multi-valued characteristics for the
subsystems. The purpose of our paper is to show that such an
extension is indeed possible. In our main result, we  prove that a
negative feedback {}interconnection of monotone i/o subsystems
with multi-valued characteristics is itself multi-stable, provided
that all the solutions of a particular discrete-time inclusion
(which is typically of much lower dimension than the subsystems)
converge.{}

Our work provides a significant extension of the Angeli-Sontag
monotone control systems theory \cite{AS03} because \cite{AS03}
requires singleton-valued characteristics and therefore globally
asymptotically stable equilibria. For other approaches to proving
multi-stability, see \cite{AS04} (where {\it positive} feedback
interconnections  of monotone i/o subsystems are considered and
the trajectories converge for {\em almost all} initial values) and
\cite{R01} (which is based on density functions and also concludes
convergence for almost all initial values).  This earlier work
does not include ours because for example (a) our results provide
global stabilization from all initial values, (b) we do not
require any regularity such as singleton-valuedness,
differentiability, or non-degeneracy for the i/s characteristics,
and (c) our results are intrinsic in the sense that we make no use
of Lyapunov or density functions.

{}{}This note is organized as follows.
In Section \ref{sec2}, we provide the necessary definitions and
background for monotone control systems, multi-valued
characteristics, weakly non-decreasing set-valued maps, and
asymptotically autonomous systems.  In Section \ref{sec3}, we
state our small-gain theorem and discuss its relationship to the
small-gain theorems in \cite{AS03, AS04, AS04a}.  In Section
\ref{sec4}, we prove our theorem and we illustrate our theorem in
Section \ref{sec5}.  We close in Section \ref{sec6} with some
suggestions for future research.

\section{Background and Motivation}
\label{sec2}
\subsection{Monotonicity and Characteristics}
\label{sec2.1} We next provide the relevant definitions for
monotone control systems and input-state characteristics.  While
our monotonicity definitions follow \cite{AS03}, our treatment of
characteristics is novel because we allow discontinuous
multi-valued characteristics and unstable equilibria. Our general
setting is that of an input-output (i/o) system
\begin{equation}\label{dyn} \dot x=f(x,u),\; \; y=h(x),\; \;
x\in {\cal X},\; \; u\in {\cal U},\; \; y\in {\cal Y}
\end{equation}
where ${\cal X}\subseteq \R^n$ is the closure of its interior and
partially ordered, ${\cal U}$ and ${\cal Y}$ are subsets of
partially ordered Euclidean spaces ${\cal B}_{\cal U}$ and ${\cal
B}_{\cal Y}$ respectively, and $f$ and $h$ are locally Lipschitz
on some open set $X$ containing ${\cal X}$. We refer to ${\cal X}$
as the {\em state space} of \rref{dyn}, ${\cal U}$ as its {\em
input space}, and ${\cal Y}$ as its {\em output space}. In
general, ${\cal X}$ will not be a linear space, since for example
we often take ${\cal X}=\R^n_{\ge 0}:=\{x\in \R^n: x_i\ge 0\,
\forall i\}$. We use $\preceq$ to denote the partial orders on all
our spaces, bearing in mind that the partial orders on our various
spaces could differ.

The set of {\em control functions} (also called {\em inputs})  for
\rref{dyn}, which we denote by ${\cal U}_\infty$, consists of all
locally essentially bounded Lebesgue measurable functions
$\uc:\R\to {\cal U}$, and we let $t\mapsto \phi(t, x_o,\uc)$
denote the trajectory of \rref{dyn} for any given initial value
$x_o\in {\cal X}$ and $\uc\in {\cal U}_\infty$.  We always assume
our dynamics $f$ are {\em forward complete} and {\em ${\cal
X}$-invariant}, which means that $\phi(\cdot, x_o,\uc)$ is defined
on $[0,\infty)$ and valued in ${\cal X}$ for all $x_o\in {\cal X}$
and $\uc\in {\cal U}_\infty$. Since we will be considering more
than one dynamic, we often use sub- or superscripts to emphasize
the state space variable or dynamic, so for example $\phi^f$ is
the flow map for the dynamic $f$ and ${\cal Y}_z$ is the output
space for an i/o system with state variable $z$.

We always assume that our  partial orders $\preceq$ are induced by
 distinguished closed nonempty sets $K$ (called {\em ordering
 cones})
 and we sometimes write
 $K_{\cal U}$ to indicate the cone inducing the partial order on
 the input space ${\cal U}$ and similarly for the other partial orders.
We always assume $K$ is a pointed convex cone, meaning,
 \[aK\subseteq K\; \;  \forall a\ge 0,\; \; \; \;
K+K\subseteq K,\; \; \; \;  K\cap (-K)=\{0\}.\] When we say that a
cone $K$ induces a
 partial order $\preceq$, we mean the following: $x\preceq y$ if
 and only if $y-x\in K$.  This induces a partial order on the set of control functions
 ${\cal U}_\infty$ as follows:  $\uc\preceq\vc$ if and only if
 $\uc(t)\preceq \vc(t)$ for Lebesgue almost all (a.a) $t\ge 0$.  A
 function $g$ mapping a partially ordered space into another
 partially ordered space is called {\em monotone} provided: $x\preceq y$
 implies $g(x)\preceq g(y)$.  We say that (\ref{dyn}) is {\em single-input single-output (SISO)} provided
${\cal B}_{\cal U}={\cal B}_{\cal Y}=\R$, taken with the usual
order, i.e., the order induced by the cone $K=[0,\infty)$.

\begin{definition}
\label{monodef} We say that \rref{dyn} is {\em monotone} provided
$h$ is monotone and
\[
(p\preceq q\; {\rm and}\;  \uc\preceq\vc)\; \; \implies\; \;
(\phi(t,p,\uc)\preceq  \phi(t,q,\vc)\; \forall t\ge 0)
\]
holds for all $p,q\in {\cal X}$ and $\uc,\vc\in {\cal
U}_\infty$.\end{definition}

We let ${\rm Equil}(f)$ denote the set of all equilibrium pairs
for our dynamic $f$,  namely, the set of all input-state pairs
$(\bar u,\bar x)$ such that $f(\bar x,\bar u)=0$. For each $(\bar
u,\bar x)\in {\rm Equil}(f)$, we let ${\cal D}^f(\bar u,\bar x)$
denote the {\em domain of attraction} of $\dot x=f(x,\bar u)$ to
$\bar x$, namely, the set of all $p\in {\cal X}$ for which
$\phi(t,p,\bar u)\to \bar x$ as $t\to +\infty$, where $\phi$ is
the flow map for $f$.  Since we are not assuming our equilibria
are stable, the sets ${\cal D}^f(\bar u,\bar x)$ are not
necessarily open and could even be singletons; see below for an
example where ${\cal D}^f(\bar u,\bar x)$ is not open. Given
$(\bar u, \bar x)\in {\rm Equil}(f)$, we say that $f$ is {\em
static Lyapunov stable at $(\bar u, \bar x)$} provided the
following condition holds for all $\eps>0$:  There exists
$\delta=\delta(\bar u,\bar x,\eps)>0$ such that for all
$x_o\in{\cal D}^f(\bar u,\bar x)\cap {\cal B}_\delta(\bar x)(=$
radius $\delta$ open ball centered at $\bar x$), we have
$|\phi(t,x_o,\bar u)-\bar x|\le \eps$ for all $t\ge 0$.

Recall the following notions from \cite{T92}, in which we let
$f^{\bar u}$ denote the constant input system $f(\cdot, \bar u)$
for each $\bar u\in {\cal U}$.  Given $\bar u\in{\cal U}$, we say
that two nonempty (but not necessarily distinct) sets
$M_1,M_2\subseteq {\cal X}$ are {\em $f^{\bar u}$-chained}
provided there exists a value $y\in {\cal X}\setminus (M_1\cup
M_2)$ and a trajectory $x:\R\to {\cal X}$ for $f^{\bar u}$
satisfying $x(0)=y$ whose {\em $\alpha$-limit set}
$\alpha(x):=\bigcap\{\overline{x((-\infty,-t])}:t\ge 0\}$ lies in
$M_1$ and whose {\em $\omega$-limit set}
$\omega(x):=\bigcap\{\overline{x([t,+\infty))}:t\ge 0\}$ lies in
$M_2$.  We {}say that a finite collection of nonempty sets $M_1,
M_2,\ldots, M_r\subseteq {\cal X}$ is {\em $f^{\bar u}$-cyclically
chained} provided the following holds:  If $r=1$, then $M_1$ is
$f^{\bar u}$-chained to itself; and if $r>1$, then $M_i$ is
$f^{\bar u}$-chained to $M_{i+1}$ for $i=1,2,\ldots, r-1$ and
$M_r$ is $f^{\bar u}$-chained to $M_1$. In this case, we call
$\{M_i\}$ an {\em $f^{\bar u}$-cycle}. An {\em $f^{\bar
u}$-equilibrium} is defined to be any point $\bar x\in {\cal X}$
 such that $f(\bar x,\bar u)=0$.  A set  $M\subseteq {\cal X}$ is called
{\em $f^{\bar u}$-invariant} provided the flow map $\phi$ for
{}$f$ satisfies $M=\{\phi(t,x,\bar u):  t\ge 0, x\in M\}$. {}A
compact $f^{\bar u}$-invariant set $M\subseteq {\cal X}$ is called
 {\em
$f^{\bar u}$-isolated compact invariant} provided there exists an
open set ${\cal U}\subseteq {\cal X}$ such that there is no
compact $f^{\bar u}$-invariant subset $\tilde M\subseteq {\cal X}$
satisfying $M\subseteq\tilde M\subseteq {\cal U}$ except $M$. {}We
use the symbol $\rightrightarrows$ to denote a {\em set-valued
map} {}{}(also called a {\em multifunction}), e.g., $F:{\cal
Z}_1\rightrightarrows{\cal Z}_2$ means that $F$ assigns each
$p\in {\cal Z}_1$ a nonempty set $F(p)\subseteq {\cal Z}_2$.

\begin{definition}\label{is}
We say that (\ref{dyn}) is endowed with a {\em static input-state
(i/s) characteristic} $k_x:{\cal U}\rightrightarrows {\cal X}$
provided:
\begin{enumerate}
\item ${\rm Graph}(k_x)={\rm Equil}(f)$; \item $\cup\{{\cal
D}^f(\bar u,\bar x): \bar x\in k_x(\bar u)\}={\cal X}$ for all
$\bar u\in {\cal U}$; \item $f$ is static Lyapunov stable at each
 $(\bar u,\bar x)\in{\rm
Equil}(f)$; and \item For each $\bar u\in {\cal U}$,  $k_x(\bar
u)$ consists of $f^{\bar u}$-isolated compact invariant $f^{\bar
u}$-equilibria and contains no $f^{\bar u}$-cycles.\end{enumerate}
In this case, we also call $k_y:=h\circ k_x$ an {\em input-output
(i/o) characteristic} for (\ref{dyn}).
\end{definition}

This definition reduces to the usual singleton-valued i/s
characteristic definition in \cite{AS03} when ${\rm
Card}\{k_x(\bar u)\}= 1$ for all $\bar u\in {\cal U}$.  We will
not use the static Lyapunov stability  property in the proof of
our small-gain theorem {\em per se},
but we still include it to make our definition of i/s
characteristics include the singleton-valued characteristic
definition in \cite{AS03}. Condition 3 in our definition is not
implied by the other conditions in the definition, even if $f$ has
no controls, since it is well-known that $f$ could admit an
unstable globally attractive equilibrium; see for example
\cite[pp. 191-4]{H67}. Condition 2 in the definition says for each
$\bar u\in {\cal U}$ and each initial state, the corresponding
$f^{\bar u}$-trajectory asymptotically approaches some state $\bar
x\in k_x(\bar u)$ (where $\bar x$ can in principle depend on the
initial state of the trajectory).  The stipulation in the static
Lyapunov stability  definition that $x_o\in{\cal D}^f(\bar u,\bar
x)\cap {\cal B}_\delta(\bar x)$ is motivated by the fact that our
domains of attraction  ${\cal D}^f(\bar u,\bar x)$ may or may not
be open, even if there are no controls. Condition $4$ is needed to
apply the theory of asymptotically autonomous systems;  see
Section \ref{sec2.3} for the relevant definitions and details.



{\em Remark:\ }Condition $4$, and in particular the ``no cycles''
part, may be hard to check in practice, at least if the system
dimension is higher than $2$, but can often be checked using
monotonicity arguments.  Consider for instance a monotone system
${\dot x}=f(x)$ having two $f$-isolated compact invariant
equilibria $p$ and $q$ and assume that $p\prec \prec q$ (where the
latter means that $q-p$ belongs to the interior of the order cone
$K$, which is assumed to be nonempty). Then there exist
neighborhoods $N_p$ and $N_q$ of $p$ and $q$ respectively such
that $n_p\prec \prec n_q$ for all $n_p\in N_p$ and $n_q\in N_q$.
We show that $\{p,q\}$ cannot be an $f$-cycle. Suppose it was a
cycle. Then there exist points $y$ and $z$ such that
$\alpha(y)=\{p\}$, $\omega(y)=\{q\}$ and $\alpha(z)=\{q\}$,
$\omega(z)=\{p\}$. It follows in particular that there exists
$T>0$ large enough such that $n_p:=\phi(-T,y)\in N_p$ and
$n_q:=\phi(-T,z)\in N_q$. Consider the strictly ordered initial
conditions $n_p \prec \prec n_q$ for the monotone system ${\dot
x}=f(x)$. Since $\omega (n_p)=\{q\}$ and $\omega (n_q)=\{p\}$,
there exists ${\tilde T>0}$ large enough so that $\phi({\tilde
T},n_p)\in N_q$ and $\phi({\tilde T},n_q)\in N_p$ and thus
$\phi({\tilde T},n_q) \prec \prec \phi({\tilde T},n_p)$, which
contradicts monotonicity of the system. The same argument can be
used to rule out cycles containing more than two equilibria, if we
assume that the equilibria are totally ordered by $\prec \prec$
(that is, either $p\prec \prec q$ or $q\prec \prec p$ whenever $p$
and $q$ are distinct equilibria).

\subsection{Weakly Non-Decreasing Set-Valued Maps}
\label{sec2.2} A basic  property of singleton-valued i/s
characteristics $k_x$ is that they are non-decreasing in the
relevant partial orders, in the sense that the following holds for
all $u,v\in {\cal U}_x$:
 $u\preceq v$ implies $k_x(u)\preceq k_x(v)$;  see \cite{AS03} for the
 elementary proof.  It is therefore natural to inquire about
 whether set-valued i/s characteristics posses some analogous
 (but more general) order-preserving property.  This motivates the
 following definition and lemma:
\begin{definition}
\label{wnddef} Let ${\cal Z}_1$ and ${\cal Z}_2$ be  partially
ordered Euclidean spaces and $F:{\cal Z}_1\rightrightarrows {\cal
Z}_2$ be any set-valued map.  We say that $F$ is {\em weakly
non-decreasing} provided the following holds for all $p,q\in {\cal
Z}_1$ such that $p\preceq q$:  For each $k_p\in F(p)$ and $k_q\in
F(q)$, there exist $r_p\in F(p)$ and $r_q\in F(q)$ such that
$r_p\preceq k_q$ and $k_p\preceq r_q$.\end{definition}
\begin{lemma}
\label{wndlem} If $k_x$ is an i/s characteristic for \rref{dyn}
and \rref{dyn} is monotone, then $k_x$ is weakly
non-decreasing.\end{lemma}
\begin{proof}
Let $p,q\in {\cal U}_x$ be such that $p\preceq q$,  let $k_p\in
k_x(p)$ and $k_q\in k_x(q)$, and let $\phi$ denote the flow map of
$f$. The corresponding trajectories for the constant inputs
satisfy $\phi(t,k_q,p)\preceq\phi(t,k_q,q)= k_q$ for all $t\ge 0$,
and $\phi(t,k_q,p)\to r_p$ for some $r_p\in k_x(p)$ as $t\to
+\infty$, so $r_p\preceq k_q$ follows because ordering cones are
closed. The other order inequality is proved similarly.
\end{proof}

Definition \ref{wnddef} reduces to non-decreasingness in the
relevant orders when $F$ is singleton-valued. We are especially
interested in  solution sequences $w_k$ satisfying discrete
set-valued inclusions $w_{k+1}\in F(w_k)$ for all $k\in \N$ where
$F$ is weakly non-decreasing. To further motivate our study of
weakly non-decreasing multifunctions, let us first assume that
$F:[0,1]\rightarrow [0,1]$ is a singleton-valued and
non-decreasing map in the usual orders (that is, $F(x)\leq F(y)$
when $x\leq y$). Then it is obvious that every solution of $
x_{k+1}=F(x_k)$ converges. Indeed, either $x_0\leq F(x_0)$ and
then $x_0\leq F(x_0)\leq F^2(x_0)\leq \dots \leq F^k(x_0)$ for all
$k\in \N$, so the sequence $\{F^k(x_0)\}$ must converge since it
is bounded above by $1$; or else $F(x_0)\leq x_0$, which leads to
a non-increasing sequence $\{F^k(x_0)\}$.  {}That converges as
well since it is bounded below by $0$. On the other hand, this
simple dynamical behavior will not occur in general for  {\it
multi-valued}, weakly non-decreasing maps.

\begin{figure}[t]
\centering
\includegraphics[width=2.5in, height=2.3in]{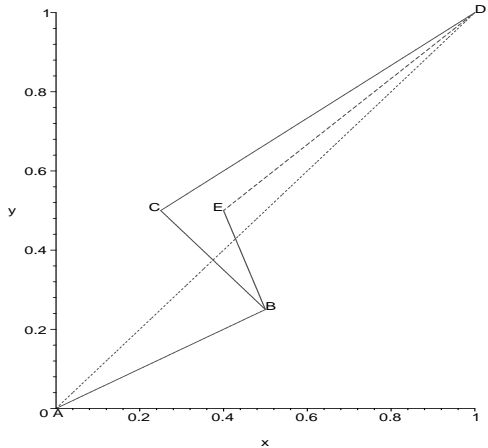}
\caption{The inverted Zorro map $F$ (ABCD) and its perturbation
$F_\epsilon$ with $\epsilon=1.5$ (ABED) from Section
\ref{sec2.2}.} \label{Zorro}
\end{figure}

To see why, consider the following simple example. Assume that
$F:[0,1]\rightrightarrows [0,1]$ is a multi-valued map whose graph
consists of the union of three straight line segments: one
connecting $A=(0,0)$ with $B=(1/2,1/4)$, a second connecting $B$
to $C=(1/4,1/2)$ (of slope $-1$), and a third connecting $C$ with
$D=(1,1)$. This ``inverted Zorro  map'' is illustrated in Figure
$\ref{Zorro}$ below and is weakly non-decreasing in the usual
orders. Then the inclusion $ x_{k+1}\in F(x_k) $ has periodic
points of period $2$. For instance, the periodic sequence
$\{1/2,1/4,1/2,1/4,...\}$ is a solution of the inclusion. In fact,
to every initial condition $x_0\in [1/4,1/2]$ corresponds a
periodic sequence of period $2$ satisfying the inclusion, namely
$\{x_0,3/4-x_0,x_0,3/4-x_0,...\}$ (since $3/4-x\in F(x)$ for all
$x\in [1/4, 1/2]$).

These periodic sequences are caused by the fact that the slope of
the middle line segment of the graph of $F$ is $-1$. Any slight
decrease of this slope will destroy the periodic points and leads
to solutions that  converge to one of the fixed points. For
example, for arbitrary $\epsilon>0$ we can define $F_{\epsilon}$
as the map whose graph consists of three straight line segments
connecting $A$ to $B$, $B$ to
$E=((1+2\epsilon)/(4+4\epsilon),1/2)$ (so the slope of
this line segment is $-1-\epsilon$), and $E$ to $D$. Then every
solution of the inclusion $x_{k+1}\in F_{\epsilon}(x_k)$ will
converge to one of the three fixed points of $F$.  In fact, each
solution sequence of this inclusion converges to either $0$ or
$1$, except for the constant sequence at the middle fixed point
$\tilde x=(3+2\epsilon)/(4(2+\epsilon))$.  To see why, notice that
if
$x_o>1/2$, then $(x_k, F_\epsilon(x_k))$ remains on the segment
$\overline{E D}$, so $x_k\uparrow 1$ by the argument for the
singleton-valued case. Similarly, if
$x_o<(1+2\epsilon)/(4+4\epsilon)$, then $(x_k, F_\epsilon(x_k))$
remains on $\overline{AB}$ so $x_k\downarrow 0$ again by the
singleton-valued case; while if $x_k$ stays in
$[(1+2\epsilon)/(4+4\epsilon),1/2]$, then
$x_{k+1}=-(1+\epsilon)x_k+\frac{3}{4}+\frac{\epsilon}{2}$ for all
$k$. Then either $x_k\equiv \tilde x$, or else
$|x_{k+1}-x_k|=(1+\epsilon)^k|x_1-x_o|\to +\infty$ as
$k\to+\infty$ which is impossible.  Therefore, either $x_k$ stays
at $\tilde x$,
 or else $x_k$ exits $[(1+2\epsilon)/(4+4\epsilon),
1/2]$ and then converges to either $0$ or $1$, as claimed.

\subsection{Asymptotically Autonomous Systems}
\label{sec2.3}

We will be especially interested in dynamics for which the
asymptotic behavior under constant inputs is known. We will then
obtain information about the trajectories for not-necessarily
constant inputs using the theory of asymptotically autonomous
systems. Before turning to this theory,
 first recall the following ``Converging-Input
Converging-State'' (CICS) Property.  This property was shown in
\cite{S03} and was used in \cite{AS03} to study the stability of
interconnected monotone systems.  We use the CICS property at the
very end of the proof of our main result (on p.\pageref{CICSuse}).
\begin{lemma}
\label{cics} Let $\bar u\in {\cal U}$, and let $\bar x$ be an
asymptotically stable equilibrium  point for $f^{\bar u}$.  Let
$\cal K$ be a compact subset of ${\cal D}^f(\bar u,\bar x)$. If
$x:[0,\infty)\to {\cal X}$ is a $\cal K$-recurrent trajectory of
$f$ for some continuous input $u:[0,\infty)\to {\cal U}$, and if
$u(t)\to \bar u$ as $t\to +\infty$, then $x(t)\to \bar x$ as $t\to
+\infty$.\end{lemma} Here $\cal K$-recurrent means for each $T>0$,
there exists $t>T$ such that $x(t)\in {\cal K}$.
One of the requirements of  asymptotic stability of $\bar x$ (in
addition to the convergence condition) is
the following stability property:  For each $\eps>0$, there exists $\delta>0$ such
that $|\phi(t,\xi,\bar u)-\bar x|\le \eps$ for all $\xi\in {\cal
B}_\delta(\bar x)$ and $t\ge 0$. The proof of the CICS property in
\cite{S03} uses  the fact that ${\cal D}^f(\bar u,\bar x)$ is
open, which follows from the assumption that $\bar x$ is a {\em
stable} equilibrium.

However, in our more general setting where the i/s characteristics
are multi-valued, the domains of attraction will not necessarily
be open, so the CICS property does not apply.  Instead, we prove
our {}result using the  theory of asymptotically autonomous
systems developed by Thieme in \cite{T92}.
  To this end, we first note that Condition 2 from our definition
  of i/s characteristics implies the
following equilibrium condition (EC)  from \cite{T92}:
\begin{itemize}\item[]\begin{itemize}
\item[(EC)] For each $\bar u\in {\cal U}$, the $\omega$-limit set
of any pre-compact $f^{\bar u}$-trajectory on $[0,\infty)$
consists of an $f^{\bar u}$-equilibrium.
\end{itemize}\end{itemize}

By an {\em asymptotically autonomous system}, we mean a system
$\dot x=H(t,x)$ that admits a second dynamic $\dot x=\bar H(x)$
(called a {\em limiting dynamic}) such that $H(t,x)\to \bar H(x)$
as $t\to +\infty$ locally uniformly in $x$.  For example, if $u\in
{\cal U}_\infty$ is continuous and $\bar u\in {\cal U}$ is such
that $u(t)\to \bar u$ as $t\to +\infty$, then for our locally
Lipschitz dynamic $f$, we know $\dot x=H(t,x):=f(x,u(t))$ is
asymptotically autonomous with limiting dynamic $\dot x =\bar
H(x):=f(x,\bar u)$. Using this observation, the following is then
immediate from \cite[Corollary 4.3]{T92}  and our i/s
characteristic definition:


\begin{lemma}
\label{T92Lemma} Assume (\ref{dyn}) is endowed with an i/s
characteristic.  Let $\bar u\in {\cal U}$ and $u:[0,\infty)\to
{\cal U}$ be any locally Lipschitz function for which $u(t)\to
\bar u$ as $t\to +\infty$.  Let $x:[0,\infty)\to{\cal X}$ be any
bounded trajectory for (\ref{dyn}) and this input $u(t)$. Then
$x(t)$ converges towards an $f^{\bar u}$-equilibrium    as $t\to
+\infty$.\end{lemma}

If one drops the  ``no cycle'' part of condition $4$ in Definition
$\ref{is}$, then the conclusion of the above Lemma does not
necessarily hold; see \cite{T92} for an example.

\section{Statement and Discussion of Small-Gain Theorem}
\label{sec3}
 We turn next to our small-gain theorem, which
generalizes
 \cite[Theorem 3]{AS03}.  The main novelty of
our result lies in its applicability to cases where one of the
interconnected systems has a multi-valued i/s characteristic, but
see Remark \ref{svis}  below for a further extension for  cases
where {\em both} subsystems  have multi-valued i/s
characteristics.
 In
what follows,  an {\em equilibrium of a discrete inclusion}
$w_{k+1}\in F(w_k)$ is defined to be any value $\bar w$ such that
$\bar w\in F(\bar w)$; the set of all equilibria for this
inclusion is denoted by ${\cal E}(F)$. A multi-function $F$ is
called {\em locally bounded} provided it maps bounded sets into
bounded sets. We say that a continuous time dynamics $F$ has a
{\em pointwise globally attractive set} $S$ provided each
maximal trajectory $\zeta(t)$ for $F$ asymptotically approaches
some point in $S$ (which could in principle depend on the specific
trajectory) as $t\to+\infty$.

\begin{theorem}
\label{sgt} Consider the following interconnection of two SISO
dynamic systems:
\begin{equation}
\label{intercon}
\begin{array}{ll}
\dot x=f_x(x,w), & y=h_x(x)\\
\dot z=f_z(z,y), & w=h_z(z)
\end{array}
\end{equation}
with ${\cal U}_x={\cal Y}_z$ and ${\cal U}_z={\cal Y}_x$.  Assume
the following:
\begin{enumerate}
\item The first system is monotone when its input $w$ and output
$y$ are ordered by the ``standard order'' induced by the positive
real semi-axis. \item The second system is monotone when its input
$y$ is ordered by the standard order and its output $w$ is ordered
by the opposite order (induced by the negative real semi-axis).
\item The respective static i/s characteristics $k_x$ and $k_z$
exist with $k_x$ singleton-valued and $k_z$ locally bounded.
 \item
Each trajectory of (\ref{intercon}) is bounded; and each solution
sequence $\{v_k\}$ of $v_{k+1}\in (k_y\circ k_w)(v_k)$ converges.
\end{enumerate} Then
(\ref{intercon}) has the pointwise globally attractive set
$\cup\{\{k_x(\bar w)\}\times (k_z\circ k_y)(\bar w): \bar w\in
{\cal E}(k_w\circ k_y)\}$.\end{theorem} In this setting,
$k_y=h_x\circ k_x$ and $k_w=h_z\circ k_z$.

Our theorem differs from  the small-gain theorem \cite[Theorem
3]{AS03} mainly in that (a)  we replaced the single valuedness of
$k_z$ with local boundedness of $k_z$, (b) we replaced the
discrete system $w_{k+1}= (k_w\circ k_y)(w_k)$ from \cite{AS03}
with a discrete inclusion, and (c)  we conclude that
(\ref{intercon}) is attracted to a set of equilibrium points
rather than a single point as in \cite{AS03}. Moreover, in
contrast to \cite{AS04}, our theorem gives {\em global}
convergence of the interconnection from all initial values.

\begin{remark}
\label{assum4} \rm  Assumption 4 of our theorem is equivalent to
the following: {\em $4'$.  Each trajectory of (\ref{intercon}) is
bounded; and $\{k_y(w_k)\}$ converges for each solution sequence
$\{w_k\}$ of $w_{k+1}\in (k_w\circ k_y)(w_k)$.}  In fact, if
Assumption 4 holds and $w_k$ is any solution of $w_{k+1}\in
(k_w\circ k_y)(w_k)$, then $k_y(w_k)$ converges because
$v_k=k_y(w_k)$ is a solution sequence for $v_{k+1}\in (k_y\circ
k_w)(v_k)$.  Conversely, if Assumption $4'$ holds, and if $v_k$ is
any solution sequence of $v_{k+1}\in (k_y\circ k_w)(v_k)$, then we
can inductively find a new sequence $r_k$ such that $v_{k+1}\equiv
k_y(r_k)$ and $r_{k+1}\in(k_w\circ k_y)(r_k)$ for all $k$, so
$v_k$ converges.  On the other hand, it could be that Assumption 4
holds but that there exists a divergent sequence $w_k$ for
$w_{k+1}\in (k_w\circ k_y)(w_k)$.  See Remark \ref{multiequil} for
an example where this occurs.  However, if the trajectories of
(\ref{intercon}) are bounded, and if each solution of
$w_{k+1}\in(k_w\circ k_y)(w_k)$ converges, then Assumption $4'$
(or equivalently Assumption $4$) holds  because $k_y$ is
continuous (by the arguments from \cite[Proposition V.5]{AS03} and
our assumption that $k_x$ is singleton valued).
\end{remark}

\section{Proof of Small-Gain Theorem}
\label{sec4}

  The following key lemma  generalizes
\cite[Proposition V.8]{AS03} to systems with multi-valued
characteristics. In it, we set  $u_{\inf}:=\liminf_{t\to
+\infty}u(t)$ and $u_{\sup}:=\limsup_{t\to +\infty}u(t)$ for any
continuous scalar function $u$ on $[0,\infty)$.

\begin{lemma}\label{keylemma}Under the hypotheses of Theorem
\ref{sgt}, if $(x(t),z(t))$ is any trajectory of (\ref{intercon})
and  $\zeta\in \omega(z)$, then there exist $k_{\scriptscriptstyle
-}\in k_z(y_{\inf})$ and $k_{\scriptscriptstyle +}\in
k_z(y_{\sup})$ such that $k_{\scriptscriptstyle -}\preceq
\zeta\preceq k_{\scriptscriptstyle +}$.
\end{lemma}
\begin{proof}
We only prove the existence of $k_-$ since the proof of the
existence of $k_+$ is similar. Set $\mu=y_{\inf}$ and let $\xi$ be
the initial value for $z(t)$. Let $t_j\to +\infty$ and $\mu_j\to
\mu$ be sequences such that  $\mu_j\in {\cal U}_z$ and $y(t)\ge
\mu_j$ for all $t\ge t_j$ and all $j$.   We have the following for
all $t\ge t_j$ and $j\in \N$:
\begin{equation}\label{old}
{}z(t)\; \; =\; \; \phi(t,\xi,y)\; \; =\; \; \phi(t-t_j,
\phi(t_j,\xi,y), y(\cdot+t_j))\; \; \succeq\; \;  \phi(t-t_j,
\phi(t_j,\xi,y), \mu_j),
\end{equation}
where $\phi$ is the flow map for $f_z$ and the last order
inequality follows from the monotonicity of the $z$-subsystem.
 Therefore, if $z(s_l)\to \zeta$ for
some sequence $s_l\to+\infty$, then we can set $t=s_l$ in
(\ref{old}) and use the closedness of order cones to find values
$v_j\in k_z(\mu_j)$ such that
\begin{equation}\label{bif}
\zeta\succeq \lim_{l\to\infty} \phi(s_l-t_j, \phi(t_j,
\xi,y),\mu_j)=v_j\; \; \; \forall j\in \N.
\end{equation}
Since $k_z$ is assumed to be locally bounded and has a closed
graph (by the continuity of the dynamic $f_z$ in all arguments),
we can find $k_-\in k_z(\mu)$ such that $\zeta\succeq v_j\to k_-$,
possibly by passing to a subsequence without relabelling.  This
proves the desired inequality.
\end{proof}

Returning to the proof of our small-gain theorem, notice that
since the output $w$ is ordered by the negative real semi-axis,
and since $k_z$ is weakly non-decreasing (by Lemma \ref{wndlem}),
it follows that
\begin{equation}
\label{antimon} \max_{k_p\in k_w(p)}\min_{k_q\in
k_w(q)}(k_p-k_q)(p-q)\; \le \;  0\; \; \; \; \forall p,q\in {\cal
U}_z.
\end{equation}
In other words, for each $p,q\in {\cal U}_z$ and $k_p\in k_w(p)$,
we can find $k_q\in k_w(q)$, such that $k_p-k_q$ and $p-q$ have
opposite signs. Also, $k_y$ is continuous and non-decreasing, as
shown in \cite[Propositions V.5 and V.8]{AS03} and Lemma
\ref{wndlem}. Choose any initial value $\xi$ for the
interconnection \rref{intercon}, and let $(x(t),z(t))$ denote the
corresponding trajectory for (\ref{intercon}) starting at $\xi$.
This trajectory is defined on $[0,\infty)$ since we are assuming
our trajectories are bounded. Set $w_{\scriptscriptstyle
+}=w_{\sup}$, $w_{\scriptscriptstyle -}=w_{\inf}$, and similarly
define $y_\pm$. Let $z_+$ (resp., $z_-$) $\in \omega(z)$ be such
that $w_-=h_z(z_+)$ (resp., $w_+=h_z(z_-)$). These limits exist
because $h_z$ is continuous and $z(t)$ is bounded in the closed
set ${\cal X}_z$. By Lemma \ref{keylemma}, we can find $k_+\in
k_z(y_+)$ and $k_-\in k_z(y_-)$ such that $k_-\preceq z_-$ and
$z_+\preceq k_+$. Setting $r^{\ssm (0)}_+=h_z(k_+)$ and $r^{\ssm
(0)}_-=h_z(k_-)$ and recalling that $w$ reverses order gives
\begin{equation} \label{key1} k_w(y_{\scriptscriptstyle +})\ni r^{\scriptscriptstyle (o)}_{\scriptscriptstyle +}\le w_{\scriptscriptstyle -}\le
w_{\scriptscriptstyle +}\le r^{\scriptscriptstyle
(o)}_{\scriptscriptstyle -}\in k_w(y_{\scriptscriptstyle
-}).\end{equation}  Since we are assuming $k_x$ is
singleton-valued, the proof of \cite[Theorem 3]{AS03} gives
\begin{equation}
\label{key2} k_y(w_{\scriptscriptstyle -})\; \le \;
y_{\scriptscriptstyle -}\; \le \;  y_{\scriptscriptstyle +}\; \le
\; k_y(w_{\scriptscriptstyle +}).
\end{equation}
Combining (\ref{key1}) and (\ref{key2}) and recalling that $k_y$
is non-decreasing gives
\begin{eqnarray}
\label{series} (k_y\circ k_w)(y_+)\; \ni\;  k_y(r^o_+)\; =:\;
s^{\scriptscriptstyle (1)}_+&\; \le \; & k_y(w_-)\; \le \; y_-\\
&\le&   y_+\;  \le \; k_y(w_+)\; \le \; s^{\scriptscriptstyle
(1)}_-\; :=\; k_y(r^{\scriptscriptstyle (o)}_-)\; \in\;  (k_y\circ
k_w)(y_-).\nonumber
\end{eqnarray}
In summary,
\begin{equation}
\label{key3} (k_y\circ k_w)(y_+)\ni s^{\scriptscriptstyle (1)}_+\;
\le \; y_-\; \le  \; y_+\; \le \; s^{\scriptscriptstyle (1)}_- \in
(k_y\circ k_w)(y_-). \end{equation} Since $y_+\; \le \;
s^{\scriptscriptstyle (1)}_-$ and $r^{\scriptscriptstyle (0)}_+\in
k_w(y_+)$, we can use (\ref{antimon}) to find
$r^{\scriptscriptstyle (1)}_+\in k_w(s^{\scriptscriptstyle
(1)}_-)\subseteq k_w(k_y\circ k_w)(y_-)$ such that
$r^{\scriptscriptstyle (1)}_+\; \le \; r^{\scriptscriptstyle
(0)}_+$.  Since $k_y$ is non-decreasing, (\ref{series}) therefore
gives
\begin{equation}
\label{first} y_-\; \ge \;  k_y(r^{\ssm (0)}_+)\; \ge \;
k_y(r^{\ssm (1)}_+)=: s^{\ssm (2)}_-\in (k_y\circ k_w)^{\ssm
2}(y_-).
\end{equation}
Similarly, since  $y_-\ge s^{\scriptscriptstyle (1)}_+$ and
$r^{\scriptscriptstyle (0)}_-\in k_w(y_-)$, we can use
(\ref{antimon}) to find $r^{\scriptscriptstyle (1)}_-\in
k_w(s^{\scriptscriptstyle (1)}_+)\subseteq k_w(k_y\circ k_w)(y_+)$
such that $r^{\scriptscriptstyle (o)}_-\; \le \;
r^{\scriptscriptstyle (1)}_-$.  Hence, (\ref{series}) also gives
\begin{equation}\label{second}
y_+\; \le \;  k_y(r^{\ssm (0)}_-)\; \le \; k_y(r^{\ssm (1)}_-)=:
s^{\ssm (2)}_+\in (k_y\circ k_w)^{\ssm 2}(y_+).
\end{equation}
Combining (\ref{first}) and (\ref{second}) gives \[(k_y\circ
k_w)^{\ssm 2}(y_-)\ni s^{\ssm (2)}_-\; \le \;  y_-\; \le \;  y_+\;
\le \; s^{\ssm (2)}_+\in (k_y\circ k_w)^{\ssm 2}(y_+).\]
 Recalling (\ref{key3}) and proceeding inductively gives sequences
 $\{s^{\ssm (r)}_\pm\}$  satisfying
 the following for all
$j\in \N$:
\begin{equation}
\label{even} (k_y\circ k_w)^{\ssm 2j}(y_-)\ni s^{\ssm (2j)}_-\;
\le \; y_-\; \le \;  y_+\; \le \; s^{\ssm (2j)}_+\in (k_y\circ
k_w)^{\ssm 2j}(y_+)
\end{equation}
\begin{equation}
\label{odd} (k_y\circ k_w)^{\ssm 2j-1}(y_+)\ni s^{\ssm (2j-1)}_+\;
\le \; y_-\; \le \;  y_+\; \le \; s^{\ssm (2j-1)}_-\in (k_y\circ
k_w)^{\ssm 2j-1}(y_-).
\end{equation}
Notice that \begin{equation}\label{general}s^{\ssm (j)}_\pm\;
\in\; (k_y\circ k_w)^{j-1}(s^{\scriptscriptstyle (1)}_\pm)\; \;
\forall j\in \N.\end{equation} Therefore, Assumption 4 from our
theorem {} provides $\bar r_\pm$ such that $s^{\ssm (j)}_\pm\to
\bar r_\pm$ as $j\to +\infty$. Letting $j\to +\infty$ in
(\ref{even}) shows that $\bar r_-\le \bar r_+$.  On the other
hand, letting $j\to +\infty$ in (\ref{odd}) gives $\bar r_+\le
\bar r_-$.  Thus,
\[\bar r_+=\bar r_-=y_+=y_-=:\bar y.\]
Applying Lemma \ref{T92Lemma} to the $z$-subsystem $f=f_z$ and the
input $u(t)=y(t) \to \bar y$ shows that $z(t)\to \bar z$ for some
$\bar z\in k_z(\bar y)$. Since $h_z$ is continuous, $w(t)$
converges as well; i.e., $w_+=w_-=:\bar w$. Therefore, $\bar
w=h_z(\bar z)\in k_w(\bar y)$ and (\ref{key2}) gives $\bar
y=k_y(\bar w)$. It follows that $\bar w\in (k_w\circ k_y)(\bar
w)$, so $\bar w\in {\cal E}(k_w\circ k_y)$. Therefore, our theorem
will follow once we show that $(x(t),z(t))$ converges to some
point in $\{k_x(\bar w)\}\times (k_z\circ k_y)(\bar w)$ as $t\to
+\infty$. To this end, first note that $x(t)\to k_x(\bar w)$ as
$t\to +\infty$ as a consequence of the CICS
property\label{CICSuse} (namely Lemma \ref{cics} above) applied to
the $x$-subsystem $f=f_x$ and the input $u(t)=w(t)\to \bar w$,
because we are assuming that $k_x$ is singleton-valued. Since
$\bar z\in k_z(\bar y)=k_z(k_y(\bar w))$, this completes the proof
of the theorem.

\begin{remark}
\label{svis} \rm  One can extend our theorem to cases where $k_x$
and $k_z$ are {\em both} multi-valued. For example, our theorem
remains true if we replace its Assumption 3 by:
\begin{itemize}{}\item[{\em $3'.$}] {\em The respective i/s characteristics
$k_x$ and $k_z$ exist and are locally bounded}.\end{itemize} In
this case the conclusion of the theorem is that our
interconnection (\ref{intercon}) has the pointwise globally
attractive set $\cup\{k_x(\bar w)\times (k_z\circ k_y)(\bar w):
\bar w\in {\cal E}(k_w\circ k_y)\}$. The proof of this alternative
formulation is similar to the proof we gave above and proceeds by
a repeated application of
\begin{equation}
\label{antimon2} \min_{k_p\in k_y(p)}\max_{k_q\in
k_y(q)}(k_p-k_q)(p-q)\; \ge \;  0\; \; \; \; \forall p,q\in {\cal
U}_x.
\end{equation}
Condition (\ref{antimon2}) follows because
 $h_x$ is
monotone and
 $k_x$ is weakly
non-decreasing.  We leave the details of the proof of this more
general version of our theorem to the reader.
\end{remark}

\section{Illustration}
\label{sec5} We next illustrate our theorem using the
interconnection
\begin{equation}
\label{interconex}
\begin{array}{ll}
\dot x=-x+5+w
, & y=x\\
\dot z=-P(z)+y,  & w=\frac{1}{1+z^2}
\end{array}
\end{equation}
evolving on $[0,\infty)\times [0,\infty)$, where
$P(z)=z(2z^2-9z+12)$.  We order $x$ and $z$ by the usual cone
$[0,\infty)$. This dynamic satisfies Conditions 1-2 from Theorem
\ref{sgt}. Replacing $w$ with $\frac{1}{1+w^2}$ in
\rref{interconex} gives the planar positive feedback system
\begin{equation}\label{ex}
\begin{array}{ll}
{\dot x}=-x+5+\frac{1}{1+w^2},& y=x \\
{\dot z}=-P(z)+y,& w=z.
\end{array}\end{equation}
{}{} If we use superscripts o to label the characteristics of our
original interconnection \rref{interconex}, and if we use $k_x$
and so on to denote the characteristics of \rref{ex}, then $
k^o_x(\frac{1}{1+w^2})\equiv k_x(w)$ and $k^o_z\equiv k_z$.  Also,
if $u_{k+1}\in (k^o_w\circ k^o_y)(u_k)$ with $u_k>0$ for all $k$,
then $w_{k+1}\in (k_w\circ k_y)(w_k)$ for all $k$ when the $w_k$'s
are chosen to satisfy
\[\frac{1}{1+w^2_k}=u_k\] for all $k\in \N$. Moreover, since the
output $w$ in (\ref{interconex}) is always positive,
 $(k^o_w\circ
k^o_y)(0)\subseteq (0,\infty)$, so $u_k>0$ for all $k\ge 1$ along
all solution sequences $\{u_k\}$ of  $u_{k+1}\in (k^o_w\circ
k^o_y)(u_k)$. Therefore, if each solution sequence $\{w_k\}$ for
$w_{k+1}\in (k_w\circ k_y)(w_k)$ converges, then each solution
sequence $\{u_k\}$ for $u_{k+1}\in (k^o_w\circ k^o_y)(u_k)$
converges as well{}, which implies the required convergence of
solutions of $v_{k+1}\in (k^o_y\circ k^o_w)(v_k)$ by Remark
\ref{assum4}.
 The fact that Condition 3 will also hold for
the original interconnection \rref{interconex} will then follow
because \rref{interconex} has the same trajectories as \rref{ex}.

It therefore remains to show that \rref{ex} satisfies Condition 3
from our theorem, that all its trajectories are bounded, and that
each solution of $w_{k+1}\in (k_w\circ k_y)(w_k)$ converges. To
this end, first note that since the outputs of both subsystems in
(\ref{ex}) are also their states, i/s and i/o characteristics
coincide for (\ref{ex})-if they exist--so we can define
\[k_1=
k_x=k_y,\; \;  k_2=k_z=k_w\] wherever the characteristics exist.
The characteristic of the first subsystem in \rref{ex} is the
singleton-valued function
\[
k_1(w)=5+\frac{1}{1+w^2},\;\; w\in \R_+,
\]
while the characteristic for the second subsystem is multi-valued
and only determined implicitly as follows: $ k_2(y)=\{z\in \R :
P(z)=y \}$ for $y\in \R_+$. A bifurcation analysis of the scalar
system ${\dot z}=-P(z)+y$, treating $y\in \R_+$ as a bifurcation
parameter, shows that $k_2(y)$ is a characteristic which is
\begin{enumerate}
\item single-valued if $y\in [0,4)$ or if $y\in (5,\infty)$. \item
triple-valued if $y\in (4,5)$. \item double-valued if $y=4$ or
$5$: $k_2(4)=\{1/2,2\}\text{ and } k_2(5)=\{1,5/2\}.$
\end{enumerate}
There are two saddle-node bifurcations, one at $y=4$ and the other
at $y=5$. The four defining properties of a characteristic (see
Definition $\ref{is}$) can indeed be readily verified: For each
$y\in \R_+$, the system ${\dot z}=-P(z)+y$ has a finite number of
isolated compact equilibria and no cycles (since the system is
scalar), and every solution converges to one of the equilibria.
It is also not hard to see that $k_2$ is locally bounded. In order
to apply Theorem \ref{sgt}, we only need to verify that
{}\rref{ex} satisfies Condition $4$ of our theorem.

To check that the trajectories of \rref{ex} (or equivalently of
\rref{interconex}) are bounded, it suffices to verify the
following:  {\em Claim (G): If {}$(x(t),z(t))$ is any trajectory
of \rref{interconex} defined on some interval $[0,T]$, then there
is
a compact set $D$ depending only on $(x(0),z(0))$ (and not on $T$)
such that $(x(t),z(t))\in D$ for all $t\in [0,T]$.}
Boundedness will follow from $(G)$ by standard results for
extendability of solutions of ODE's. To prove $(G)$, first note
that the boundedness of $w$ on $[0,T]$ and the variations of
parameters formula gives \[|y(t)|=|x(t)|\le |x(0)|+6\] for all
$t\in[0,T]$. Pick $\tilde z>5/2$ such that $\tilde
z=P^{-1}(|x(0)|+6)$ which exists because $P$ is one-to-one above
$5/2$.  It follows that if {}$t\in [0,T)$ is such that
$z(t)>\tilde z$, then \[(z(t))\ge P(\tilde z)=|x(0)|+6\ge y(t),\]
so $\dot z(t)\le 0$. Therefore, $z(t)$ stays below $\tilde z$ on
$[0,T]$. Since $\tilde z$ depends only on $x(0)$, Claim (G)
follows.

 Next
consider the discrete inclusion $w_{k+1}\in \left(k_2\circ
k_1\right)(w_k)$ and notice that it reduces to a discrete equation
$ w_{k+1}=\left(k_2\circ k_1\right)(w_k) $ because $k_1(w)>5$ and
$k_2(y)$ is single-valued when $y>5$. Notice also that for all
$w_0\in \R_+$, the discrete equation gives  $w_k> 5/2$ for all
$k\geq 1$. In particular, the interval $(5/2,\infty)$ is forward
invariant for the discrete equation. Finally, since $|k'_1(w)|$ is
decreasing for $w\ge 5/2$,
 elementary calculus shows that
\[
\left| k_2'(k_1(w))k_1'(w) \right|\; \le \;
\frac{|k'_1(5/2)|}{P'(k_2\circ k_1(w))}\; \le \;
\frac{|k'_1(5/2)|}{P'(5/2)}\; =\;
\frac{5}{(1+25/4)^2}\frac{2}{9}\; <\; 1\; \; \forall w\ge 5/2,\]
so $k_2\circ k_1$ is a contraction mapping on $[5/2,\infty)$,
hence the discrete equation has a unique globally attractive fixed
point ${\bar w}$. Therefore, we know {}from Remark \ref{assum4}
that (\ref{ex}) satisfies Conditions 3-4 of our theorem, as
claimed. Since
\[
{\cal E}(k^o_w\circ k^o_y)=\left\{ \frac{1}{1+\bar w^2}: \bar w\in
{\cal E}(k_2\circ k_1)\right\},
\]
we conclude that our original interconnection \rref{interconex}
has the unique globally attractive equilibrium
\[{}{}
\left\{\left(5+\frac{1}{1+\bar w^2}, k_2\left(5+\frac{1}{1+\bar
w^2}\right)\right)\right\}.
\]
Figure \ref{graph} below illustrates this.


\begin{figure}[H]
\centering
\includegraphics[width=3.8in,height=3.5in]{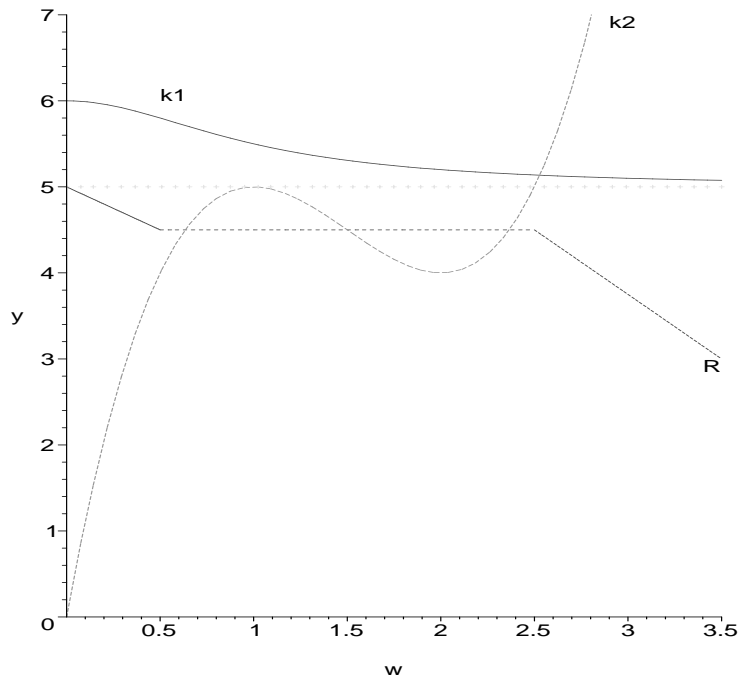}
\caption{Characteristics $k_1(w)$, $k_2(y)$ and $R(w)$ from Section
\ref{sec5}.} \label{graph}
\end{figure}

\begin{remark}
\label{multiequil} \rm In the preceding example, the inclusion
$w_{k+1}\in (k_w\circ k_y)(w_k)$ had a unique equilibrium, but our
theory applies to examples where ${\cal E}(k_w\circ k_y)$ has more
than one element as well.  One such example is constructed by
modifying the interconnection (\ref{ex}) in the following way:
replace the $x$-subsystem with $\dot x=-x+R(w)$ where $R$ consists
of the line segments in the $wy$-plane joining $(0,5)$ to
$(.5,4.5)$, $(.5, 4.5)$ to $(2.5, 4.5)$, and $(2.5, 4.5)$ to
$(3.5, 3)$.  With this change we get {}${\rm Card}\{{\cal
E}(k_w\circ k_y)\}=3$, and the conclusion of our theorem remains
true because $\{k_y(w_k)\}$ converges for each solution sequence
$\{w_k\}$ of $w_{k+1}\in (k_w\circ k_y)(w_k)$;  see Remark
\ref{assum4}.  In fact, if $w_o\in [.5,2.5]$, then $k_y(w_o)=4.5$,
so $w_k\in {\cal E}(k_w\circ k_y)$ for all $k\in \N$, which gives
$k_y(w_k)=4.5$ for all $k\in \N$.  If $w_o\in [0,.5]$, then
$k_y(w_o)\in [4.5,5]$, which gives $w_1\in [.5,2.5]$, so
$k_y(w_k)\equiv 4.5$ for all $k\ge 2$ as before.  Finally, if
$w_o>5/2$, then $k_y(w_o)\le 4.5$, so $w_1\in k_w\circ k_y(w_o)\in
[0,5/2]$, so $k_y(w_k)=4.5$ for $k\ge 3$, by the previous two
cases. On the other hand, one can find non-periodic divergent
solution sequences of $w_{k+1}\in (k_w\circ k_y)(w_k)$ when
$w_o\in {}[1/2,5/2]$. The {}detailed analysis of this more
complicated example is similar to the analysis of (\ref{ex}) and
is left to the reader.
Note that the convergence of the iterations in the preceding remark follows
because $R(w)$ is a horizontal line, at least locally where it
meets the other characteristic.
\end{remark}



\section{Conclusion}
\label{sec6} We presented a new small-gain theorem for
interconnections of  monotone i/o systems with set-valued i/s
characteristics.  This corresponds to situations where the
{}trajectory for a given constant input can converge to several
possible equilibria, depending on the initial value for the
trajectory.  {}{} A key ingredient in the proof of our small-gain
theorem is the theory of asymptotically autonomous systems, which
requires in particular that the equilibria of the subsystems in
the interconnection contain no chains.  This suggests the question
of how one might extend our theory to cases where the sets of
equilibria of the subsystems are more general, e.g., where they
contain chains or limit cycles.  Research on this question is
ongoing.


\end{document}